\documentclass[oneside,10pt,]{amsart}
\usepackage{geometry}
\usepackage{ifthen}
\usepackage{etoolbox}
\usepackage{ifxetex,ifluatex}
\usepackage{graphicx}
\PassOptionsToPackage{dvipsnames,svgnames,table}{xcolor}
\usepackage{xcolor}
\theoremstyle{plain}
\newtheorem{thmbox}{}[section]
\newtheorem{claim}[thmbox]{Claim}
\newtheorem{proposition}[thmbox]{Proposition}

\theoremstyle{definition}
\newtheorem{example}[thmbox]{Experiment}

\theoremstyle{remark}

\usepackage{tcolorbox}
\tcbuselibrary{skins}
\tcbuselibrary{breakable}
\tcbuselibrary{raster}
\tcbset{ bwminimalstyle/.style={size=minimal, boxrule=-0.3pt, frame empty,
colback=white, colbacktitle=white, coltitle=black, opacityfill=0.0} }
\tcbset{ runintitlestyle/.style={fonttitle=\blocktitlefont\upshape\bfseries, attach title to upper} }
\tcbset{ exercisespacingstyle/.style={before skip={1.5ex plus 0.5ex}} }
\tcbset{ blockspacingstyle/.style={before skip={2.0ex plus 0.5ex}} }
\tcbuselibrary{xparse}
\usetikzlibrary{calc}
\tcbuselibrary{hooks}
\newlength{\normalparskip}
\AtBeginDocument{\setlength{\normalparindent}{\parindent}}
\AtBeginDocument{\setlength{\normalparskip}{\parskip}}
\newcommand{\setparstyle}{\setlength{\parindent}{\normalparindent}\setlength{\parskip}{\normalparskip}}
\usepackage{xstring}
\ifthenelse{\boolean{xetex} \or \boolean{luatex}}{%
\ifxetex\usepackage{xltxtra}\fi
\ifluatex\usepackage{realscripts}\fi
}{
\usepackage[utf8]{inputenc}
}

\newcommand{\blocktitlefont}{\relax}

\newcommand{\xreffont}{\relax}

\ifthenelse{\boolean{xetex} \or \boolean{luatex}}{%
\usepackage{fontspec}
\IfFontExistsTF{Inconsolatazi4-Regular.otf}{}{\GenericError{}{The font "Inconsolatazi4-Regular.otf" requested by PreTeXt output processed by the  xelatex  executable is not available.  Either a file cannot be located in default locations via a filename, or a font is not known by its name as part of your system.}{Consult the PreTeXt Guide for help with LaTeX fonts, or perhaps try using  pdflatex  as a test.}{}}
\IfFontExistsTF{Inconsolatazi4-Bold.otf}{}{\GenericError{}{The font "Inconsolatazi4-Bold.otf" requested by PreTeXt output processed by the  xelatex  executable is not available.  Either a file cannot be located in default locations via a filename, or a font is not known by its name as part of your system.}{Consult the PreTeXt Guide for help with LaTeX fonts, or perhaps try using  pdflatex  as a test.}{}}
\usepackage{zi4}
\setmonofont[BoldFont=Inconsolatazi4-Bold.otf,StylisticSet={1,3}]{Inconsolatazi4-Regular.otf}
\usepackage{polyglossia}
\setmainlanguage[variant=usmax]{english}
}{%
\usepackage{lmodern}
\usepackage[T1]{fontenc}
\usepackage[varqu,varl]{inconsolata}
}
\usepackage{microtype}
\usepackage{amsmath}
\usepackage{amscd}
\usepackage{amssymb}
\allowdisplaybreaks[4]
\newcommand{\mono}[1]{\texttt{#1}}
\newcommand{\terminology}[1]{\textbf{#1}}
\numberwithin{equation}{section}
\tcbset{ imagestyle/.style={bwminimalstyle} }
\NewTColorBox{tcbimage}{mmm}{imagestyle,left skip=#1\linewidth,width=#2\linewidth}
\NewDocumentEnvironment{image}{mmmm}{\notblank{#4}{\leavevmode\nopagebreak\vspace{#4}}{}\begin{tcbimage}{#1}{#2}{#3}}{\end{tcbimage}%
}
\ifthenelse{\boolean{xetex} \or \boolean{luatex}}%
  {\tcbuselibrary{listings}}%
  {\tcbuselibrary{listingsutf8}}%
\makeatletter
\lst@CCPutMacro\lst@ProcessOther {"2D}{\lst@ttfamily{-{}}{-{}}}
\@empty\z@\@empty
\makeatother
\lstdefinelanguage{none}{identifierstyle=,commentstyle=,stringstyle=,keywordstyle=}
\ifthenelse{\boolean{xetex}}{}{%
\lstset{extendedchars=true}
\lstset{literate={ }{{~}}{1}{¡}{{\textexclamdown }}{1}{¢}{{\textcent }}{1}{£}{{\textsterling }}{1}{¤}{{\textcurrency }}{1}{¥}{{\textyen }}{1}{¦}{{\textbrokenbar }}{1}{§}{{\textsection }}{1}{¨}{{\textasciidieresis }}{1}{©}{{\textcopyright }}{1}{ª}{{\textordfeminine }}{1}{«}{{\guillemotleft }}{1}{¬}{{\textlnot }}{1}{­}{{\-}}{1}{®}{{\textregistered }}{1}{¯}{{\textasciimacron }}{1}{°}{{\textdegree }}{1}{±}{{\textpm }}{1}{²}{{\texttwosuperior }}{1}{³}{{\textthreesuperior }}{1}{´}{{\textasciiacute }}{1}{µ}{{\textmu }}{1}{¶}{{\textparagraph }}{1}{·}{{\textperiodcentered }}{1}{¸}{{\c{}}}{1}{¹}{{\textonesuperior }}{1}{º}{{\textordmasculine }}{1}{»}{{\guillemotright }}{1}{¼}{{\textonequarter }}{1}{½}{{\textonehalf }}{1}{¾}{{\textthreequarters }}{1}{¿}{{\textquestiondown }}{1}{À}{{\`{A}}}{1}{Á}{{\'{A}}}{1}{Â}{{\^{A}}}{1}{Ã}{{\~{A}}}{1}{Ä}{{\"{A}}}{1}{Å}{{\AA }}{1}{Æ}{{\AE }}{1}{Ç}{{\c{C}}}{1}{È}{{\`{E}}}{1}{É}{{\'{E}}}{1}{Ê}{{\^{E}}}{1}{Ë}{{\"{E}}}{1}{Ì}{{\`{I}}}{1}{Í}{{\'{I}}}{1}{Î}{{\^{I}}}{1}{Ï}{{\"{I}}}{1}{Ð}{{\DH }}{1}{Ñ}{{\~{N}}}{1}{Ò}{{\`{O}}}{1}{Ó}{{\'{O}}}{1}{Ô}{{\^{O}}}{1}{Õ}{{\~{O}}}{1}{Ö}{{\"{O}}}{1}{×}{{\texttimes }}{1}{Ø}{{\O }}{1}{Ù}{{\`{U}}}{1}{Ú}{{\'{U}}}{1}{Û}{{\^{U}}}{1}{Ü}{{\"{U}}}{1}{Ý}{{\'{Y}}}{1}{Þ}{{\TH }}{1}{ß}{{\ss }}{1}{à}{{\`{a}}}{1}{á}{{\'{a}}}{1}{â}{{\^{a}}}{1}{ã}{{\~{a}}}{1}{ä}{{\"{a}}}{1}{å}{{\aa }}{1}{æ}{{\ae }}{1}{ç}{{\c{c}}}{1}{è}{{\`{e}}}{1}{é}{{\'{e}}}{1}{ê}{{\^{e}}}{1}{ë}{{\"{e}}}{1}{ì}{{\`{\i}}}{1}{í}{{\'{\i}}}{1}{î}{{\^{\i}}}{1}{ï}{{\"{\i}}}{1}{ð}{{\dh }}{1}{ñ}{{\~{n}}}{1}{ò}{{\`{o}}}{1}{ó}{{\'{o}}}{1}{ô}{{\^{o}}}{1}{õ}{{\~{o}}}{1}{ö}{{\"{o}}}{1}{÷}{{\textdiv }}{1}{ø}{{\o }}{1}{ù}{{\`{u}}}{1}{ú}{{\'{u}}}{1}{û}{{\^{u}}}{1}{ü}{{\"{u}}}{1}{ý}{{\'{y}}}{1}{þ}{{\th }}{1}{ÿ}{{\"{y}}}{1}}
}
\definecolor{identifiers}{rgb}{0.375,0,0.375}
\definecolor{comments}{rgb}{0.5,0,0}
\definecolor{strings}{rgb}{0,0.5,0}
\definecolor{keywords}{rgb}{0,0,0.5}
\lstdefinestyle{programcodestyle}{identifierstyle=\color{identifiers},commentstyle=\color{comments},stringstyle=\color{strings},keywordstyle=\color{keywords}, breaklines=true, breakatwhitespace=true, columns=fixed, extendedchars=true, aboveskip=0pt, belowskip=0pt}
\lstdefinestyle{programcodenumberedstyle}{style=programcodestyle, numbers=left}
\tcbset{ programboxstyle/.style={left=3ex, right=0pt, top=0ex, bottom=0ex, middle=0pt, toptitle=0pt, bottomtitle=0pt, boxsep=0pt, 
listing only, fontupper=\small\ttfamily,
colback=white, sharp corners, boxrule=-0.3pt, leftrule=0.5pt, toprule at break=-0.3pt, bottomrule at break=-0.3pt,
breakable,
} }
\tcbset{ programboxnumberedstyle/.style={programboxstyle, left=6ex} }
\newtcblisting{program}[4]{programboxstyle, left skip=#2\linewidth, width=#3\linewidth, listing options={language=#1, style=programcodestyle}}
\newtcblisting{programnumbered}[4]{programboxnumberedstyle, left skip=#2\linewidth, width=#3\linewidth, listing options={language=#1, style=programcodenumberedstyle}}
\usepackage{enumitem}
\newlist{referencelist}{description}{4}
\setlist[referencelist]{leftmargin=!,labelwidth=!,labelsep=0ex,itemsep=1.0ex,topsep=1.0ex,partopsep=0pt,parsep=0pt}
\usepackage[hyperfootnotes=false]{hyperref}
\hypersetup{colorlinks=true,linkcolor=blue,citecolor=blue,filecolor=blue,urlcolor=blue}
\hypersetup{hypertexnames=false}
\makeatletter
\patchcmd\Hy@EveryPageBoxHook{\Hy@EveryPageAnchor}{\Hy@hypertexnamestrue\Hy@EveryPageAnchor}{}{\fail}
\makeatother
\hypersetup{pdftitle={The bimodal distribution in the derivative of unitary polynomials}}

\newtcolorbox[auto counter, number within=section]{block}{}
\newtcolorbox[auto counter, number within=section]{project-distinct}{}
\makeatletter
\newtcolorbox[auto counter, number within=tcb@cnt@block, number freestyle={\noexpand\thetcb@cnt@block(\noexpand\alph{\tcbcounter})}]{subdisplay}{}
\makeatother
\tcbset{ listingptxstyle/.style={bwminimalstyle, middle=1ex, blockspacingstyle, fontlower=\blocktitlefont} }
\newtcolorbox[use counter from=block]{listingptx}[4]{lower separated=false, before lower={{\textbf{#1~\thetcbcounter}\space#2}}, phantomlabel={#3}, unbreakable, listingptxstyle, }
\tcbset{ figureptxstyle/.style={bwminimalstyle, middle=1ex, blockspacingstyle, fontlower=\blocktitlefont} }
\newtcolorbox[use counter from=block]{figureptx}[4]{lower separated=false, before lower={{\textbf{#1~\thetcbcounter}\space#2}}, phantomlabel={#3}, unbreakable, figureptxstyle, }
\usepackage{tikz}
\tcbset{ sbsstyle/.style={raster before skip=2.0ex, raster equal height=rows, raster force size=false, raster after skip=0.7\baselineskip} }
\tcbset{ sbspanelstyle/.style={bwminimalstyle, fonttitle=\blocktitlefont, before upper app={\setparstyle}} }
\NewDocumentEnvironment{sidebyside}{mmmm}
  {\begin{tcbraster}
    [sbsstyle,raster columns=#1,
    raster left skip=#2\linewidth,raster right skip=#3\linewidth,raster column skip=#4\linewidth]}
  {\end{tcbraster}}
\NewTColorBox{sbspanel}{mO{top}}{sbspanelstyle,width=#1\linewidth,valign=#2}
\usepackage{extpfeil}

\newcommand{\RR}{\mathbb{R}}
\newcommand{\CC}{\mathbb{C}}
\newcommand{\abs}[1]{|#1|}
\newcommand{\lt}{<}

\begin{document}
\title[The bimodal distribution in the derivative of unitary polynomials]{The bimodal distribution in the derivative of unitary polynomials}
\author{David W. Farmer}
\address{American Institute of Mathematics\\
}
\email{\href{mailto:farmer@aimath.org}{\nolinkurl{farmer@aimath.org}}}
\begin{abstract}
The derivative of a polynomial with all zeros on the unit circle has the zeros of its derivative on or inside the unit circle.  It has been observed that in many cases the zeros of the derivative have a bimodal distribution: there are two smaller circles near which it is more likely to find those zeros.  We identify the likely source of the second mode.  This idea is supported with numerical examples involving the characteristic polynomials of random unitary matrices.%
\end{abstract}
\maketitle
A \terminology{unitary polynomial} is a polynomial \(p\) having all roots on the unit circle, \(\{z \in \CC \ :\ \abs{z} = 1\}\). We denote those zeros by \(e^{i\theta_1},\ldots,e^{i\theta_N}\), where \(0\le \theta_1\le \cdots \le \theta_N \lt 2\pi\). By the Gauss-Lucas theorem, all the roots of the derivative \(p'\) lie in the unit disc,  \(\{z \in \CC \ :\ \abs{z} \le 1\}\), with the only zeros of \(p'\) on the unit circle arising from multiple zeros of~\(p\).%
\par
We consider the case of unitary polynomials with zeros \(z_1,\ldots,z_n\) randomly, but not necessarily independently, distributed on the unit circle. The zeros of the derivative \(z'_1,\ldots,z'_{n-1}\) will be random points in the unit disc: we want to understand their distribution. We will further assume that the distribution of zeros is rotationally invariant: the likelihood of zeros \(\{e^{i\theta_1},\ldots,e^{i\theta_n}\}\) is the same as the likelihood of \(\{e^{i\theta_1+\phi},\ldots,e^{i\theta_n+\phi}\}\) for any \(\phi\in\RR\). In this case it is natural to look at the radial distribution of the zeros of the derivative, that is, the distribution of~\(\abs{z'_j}\).%
\par
It is also helpful to imagine that any finite sequence of normalized neighbor spacings can occur.  See \hyperref[sec_stable]{Subsection~{\xreffont\ref{sec_stable}}}.%
\par
The zeros of \(p'\) will cluster near the unit circle, equivalently, most of the \(\abs{z'}\) will be very close to~\(1\). As described in [\hyperlink{mez}{{\xreffont 3}}, \hyperlink{dffhmp}{{\xreffont 1}}], for degree \(N\) polynomials the appropriate rescaling is to consider \(N (1 - \abs{z'})\).  It was found \hyperlink{dffhmp}{[{\xreffont 1}]} that in a wide variety of cases the distribution of \(\abs{z'}\) is bimodal: the probability density function (PDF) has two local maxima.%
\par
We perform various experiments which indicate the source of the second mode.   We focus on the case of characteristic polynomials of random unitary matrices.  That is, \(u\in U(N)\) is selected randomly with respect to Haar measure, and \(p_u\) is its characteristic polynomial. %
\subparagraph*{The three regimes.}\label{intro-6}%
We suggest that a productive way to think about the zeros of \(p'\) is to recognize three regions.  One is the annulus between the unit circle and (approximately) the circle of radius \(1-2/N\).  In that region, most of the zeros of \(p'\) arise from closely spaced zeros of~\(p\). This situation was analyzed in detail in~\hyperlink{dffhmp}{[{\xreffont 1}]}. But note that zeros of \(p'\) very close to the unit circle also arise from sufficiently long strings of consecutive zeros of \(p\) spaced closer than average: see~\hyperlink{farmer_ki}{[{\xreffont 2}]}.%
\par
The second is the annulus between (approximately) the circles of radius \(1-2/N\) and of radius \(1-4/N\).  Here some of the zeros of \(p'\) arise from zeros spaced slightly closer than average, or from clusters of such zeros.  But in that region many of the zeros of \(p'\) arise from a different source: triples of consecutive zeros, where both neighbor gaps are larger than average. See \hyperref[thesource]{Section~{\xreffont\ref{thesource}}} for a discussion.%
\par
Furthermore, in both of those regions typically the location of the zeros of \(p'\) is primarily determined by a few nearby zeros of~\(p\). That is, if there is a zero \(z'\) of \(p'_0\) in that region, and \(p_1\) is another polynomial of the same degree, with \(p_0\) and \(p_1\) having the same zeros near  \(z'\), then \(p'_1\) will typically have a zero close to \(z'\).  In other words, zeros of \(p'\) which are closer than (approximately) \(4/N\) of the unit circle, depend primarily on the local arrangement of zeros of~\(p\).%
\par
In contrast, the zeros of \(p'\) farther than (approximately) \(4/N\) from the unit circle are sensitive to the bulk of the zeros of \(p\). See \hyperref[jiggled]{Figure~{\xreffont\ref{jiggled}}} for an illustration of this phenomenon. \hyperref[combinedpdf]{Figure~{\xreffont\ref{combinedpdf}}} indicates that approximately~\(25\%\) of the zeros of \(p'\) fall into that category.%
\par
In the examples where we plot zeros of \(p'\) in the unit circle, we also show (as dotted lines) the circles of radius \(1-2/N\) and \(1-4/N\).%
\par\medskip
The existence of multiple regions with different behaviors suggests that determining the distribution of zeros of \(p'\) is truly a difficult problem.%
%
%
\typeout{************************************************}
\typeout{Section 1 Random characteristic polynomials}
\typeout{************************************************}
\section{Random characteristic polynomials}\label{pdfs}
Motivated by the connection to the Riemann \(\zeta\)-function and other L-functions, and building on previous work~[\hyperlink{mez}{{\xreffont 3}}, \hyperlink{dffhmp}{{\xreffont 1}}], our examples will primarily come from the characteristic polynomials of random unitary matrices, where the randomness is uniform with respect to Haar measure on the unitary group~\(U(n)\).%
%
%
\typeout{************************************************}
\typeout{Subsection 1.1 Some experiments to gain intuition}
\typeout{************************************************}
\subsection{Some experiments to gain intuition}\label{pdfs-3}
Generating Haar-random unitary matrices is easy in principle, but there are two likely pitfalls.  The first, explained in~\hyperlink{Mez_haar}{[{\xreffont 4}]}, is that the built-in \(QR\)-decomposition routine in many computer-algebra systems introduces a bias.  The second is that the built-in routine to find the eigenvalues of a matrix may give incorrect results if the matrix is large and the entries are not known with sufficient precision. This Mathematica code addresses both issues, although it is left to the user to ensure that sufficient precision is selected.%
\begin{listingptx}{Listing}{Mathematica code to generate Haar-random unitary matrices with a given precision.}{pdfs-3-3}{}%
\begin{program}{none}{0}{1}{0}
norm = NormalDistribution[0, 1];

randunitary[n_, prec_] := Block[{q, r, d, i, j},
   {q, r} = QRDecomposition[
      Table[RandomVariate[norm, WorkingPrecision -> prec] + 
        RandomVariate[norm, WorkingPrecision -> prec]*I,
          {i, 1, n}, {j, 1, n}]
            ];
   d = DiagonalMatrix[
         Table[dd = r[[i, i]]; dd/Abs[dd], {i, 1, n}]
       ];
   q.d]
\end{program}
\tcblower
\end{listingptx}%
Our first experiment is to compile data on the radial distribution of zeros of derivatives.%
\begin{example}[Producing a large data set.]\label{pdfs-3-5}%
Start with an empty list named \mono{zp20}. Generate \(5.4 \times 10^{6}\) Haar-random matrices in \(U(20)\).  For each matrix:%
\begin{enumerate}
\item{}form the characteristic polynomial%
\item{}take its derivative%
\item{}compute approximations to the zeros of the derivative%
\item{}append the absolute values of those zeros to the list \mono{zp20}.%
\end{enumerate}
\par\smallskip%
\noindent\textbf{\blocktitlefont Outcome}.\hypertarget{pdfs-3-5-3}{}\quad{}The list \mono{zp20} will contain approximately \(10^{8}\) real numbers \(\{r_j\}\) with \(0 \lt r_j \lt 1\). Note that those inequalities are strict, because the probability of hitting either extreme is~\(0\). \hyperref[zp20pdf1]{Figure~{\xreffont\ref{zp20pdf1}}} shows a histogram, normalized to have area~\(1\), of \(20 (1 - r_j)\).%
\end{example}
\begin{figureptx}{Figure}{Empirical distribution of \(20(1 - \abs{z'})\) for \(z'\) a zero of the derivative of the characteristic polynomial of a Haar-random matrix from \(U(20)\).}{zp20pdf1}{}%
\begin{image}{0.1}{0.8}{0.1}{}%
\includegraphics[width=\linewidth]{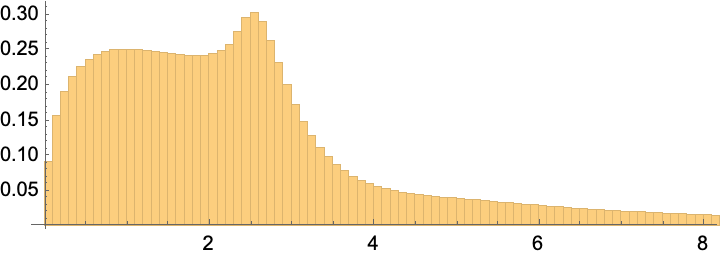}
\end{image}%
\tcblower
\end{figureptx}%
Is the second bump a transient phenomenon, or will it persist in the large \(N\) limit? \hyperref[combinedpdf]{Figure~{\xreffont\ref{combinedpdf}}} superimposes the results of the above experiment for various~\(N\).%
\begin{figureptx}{Figure}{Empirical distribution of \(N(1 - \abs{z'})\) for \(z'\) a zero of the derivative of the characteristic polynomial of a Haar-random matrix from \(U(N)\), for \(N=20\), \(40\), \(80\), \(160\), and \(320\). The dotted lines divide the \(N=320\) PDF into \(8\) equal areas.}{combinedpdf}{}%
\begin{image}{0.1}{0.8}{0.1}{}%
\includegraphics[width=\linewidth]{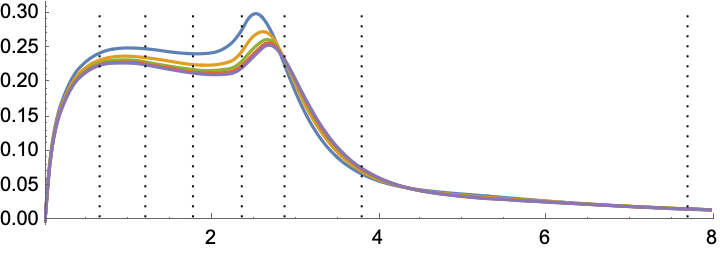}
\end{image}%
\tcblower
\end{figureptx}%
Does \hyperref[combinedpdf]{Figure~{\xreffont\ref{combinedpdf}}} support a claim that the second bump persists as~\(N\to\infty\)?%
\par
One possible argument against such interpretation is that some features, for example the value distribution of~\(p_u(1)\) for random \(u\in U(N)\), approach their limiting behavior on the scale of \(1/\sqrt{\log N}\).  For such features, most numerical experiments will show significant deviations from the limiting behavior.%
\par
We suggest that for the normalized distribution \(N(1-\abs{z'})\), the limiting distribution is approached on the scale of \(1/N^2\), and so \hyperref[combinedpdf]{Figure~{\xreffont\ref{combinedpdf}}} does provide persuasive evidence for a second bump in the limit.  There are three ingredients to our reasoning.%
\par
First is that the zeros of \(p'\) near the two peaks typically depend strongly only on the relative spacing of a few nearby zeros.  See \hyperref[sec_local]{Subsection~{\xreffont\ref{sec_local}}}.%
\par
Second is the the relative likelihood of a given cluster of zero spacings stabilizes quickly as \(N\) increases, at the rate of \(1/N^2\).  See \hyperref[sec_stable]{Subsection~{\xreffont\ref{sec_stable}}}.%
\par
Third is that the effect of the curvature of the circle on the rescaled zeros also decreases as \(1/N^2\).  See \hyperref[sec_curve]{Subsection~{\xreffont\ref{sec_curve}}}.%
%
%
\typeout{************************************************}
\typeout{Subsection 1.2 The local arrangement is most important, when close to the unit circle}
\typeout{************************************************}
\subsection{The local arrangement is most important, when close to the unit circle}\label{sec_local}
\hyperref[jiggled]{Figure~{\xreffont\ref{jiggled}}} shows the results of the following:%
\begin{example}[The sensitivity to distant zeros.]\label{p40perturbed}%
Start with an empty list named \mono{zp40perturbed}. Generate a unitary polynomial \(p_u\) for some fixed \(u\in U(40)\). Perform the following steps  \(250\) times:%
\begin{enumerate}
\item{}Generate a polynomial \(\tilde{p}\), where the zeros of \(\tilde{p}\) are partitioned into two sets:%
\begin{enumerate}
\item{}The zeros of \(p_u\) in the first quadrant%
\item{}The other zeros of \(p_u\), with arguments perturbed independently by a random amount in the interval \((-2\pi/40,2\pi/40)\).  That is, perturbed by a random amount less than the average gap between zeros.%
\end{enumerate}
\item{}Compute approximations to the zeros of \(\tilde{p}'\)%
\item{}Append those zeros to the list \mono{zp40perturbed}.%
\end{enumerate}
\par\smallskip%
\noindent\textbf{\blocktitlefont Outcome}.\hypertarget{p40perturbed-3}{}\quad{}The list \mono{zp40perturbed} will contain \(9750\) points, which are shown as small blue dots in \hyperref[jiggled]{Figure~{\xreffont\ref{jiggled}}}.  The roots of \(p_u\) are shown as larger red squares.%
\begin{figureptx}{Figure}{The outcome of one run of \hyperref[p40perturbed]{Experiment~{\xreffont\ref{p40perturbed}}}.}{jiggled}{}%
\begin{image}{0.15}{0.7}{0.15}{}%
\includegraphics[width=\linewidth]{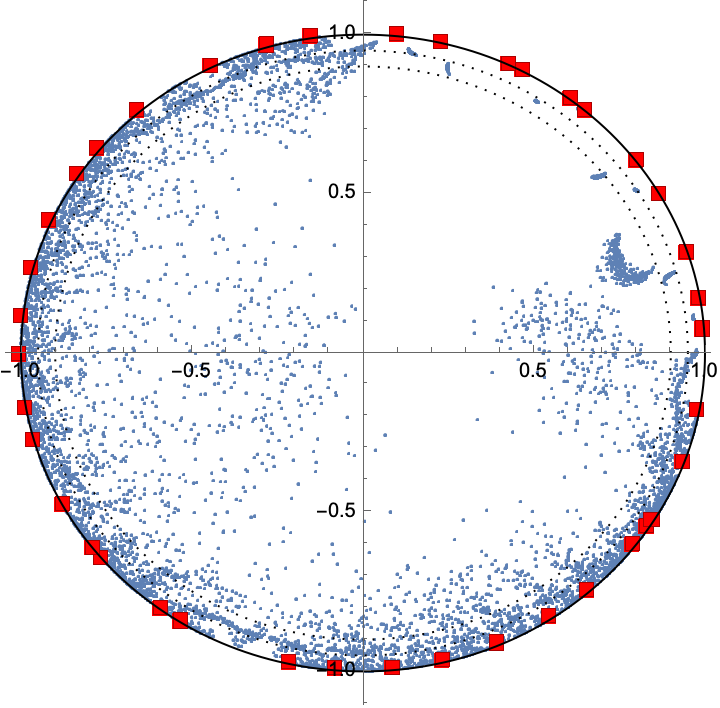}
\end{image}%
\tcblower
\end{figureptx}%
\end{example}
In \hyperref[jiggled]{Figure~{\xreffont\ref{jiggled}}} it can be seen that the zeros of the derivative which are in the middle of the first quadrant and closer than around \(4/N\) from the unit circle, depend minimally on the other zeros. Thus, only the arrangement of small numbers of consecutive zeros determines the main features of the zeros of the derivative in that region.%
%
%
\typeout{************************************************}
\typeout{Subsection 1.3 Relative frequency of neighbor spacings}
\typeout{************************************************}
\subsection{Relative frequency of neighbor spacings}\label{sec_stable}
The second ingredient is a property of Haar measure on \(U(N)\). Suppose \(\delta_1 \lt \cdots \lt \delta_k\).  We interpret those numbers as the consecutive normalized gaps between random eigenvalues. That is, we consider the likelihood that for some \(\theta_0\) we have \(\{e^{i\theta_0}, e^{i(\theta_0 + 2\pi \delta_1/N)},\ldots, e^{i(\theta_0+2\pi \delta_k/N)}\}\) as consecutive eigenvalues of a Haar random matrix in~\(U(N)\). By \terminology{likelihood} we mean%
\begin{equation}
\lim_{\epsilon \to 0} \epsilon^{-k} \operatorname{Prob}(\delta_1 \lt \tilde\lambda_1-\tilde\lambda_0 \lt \delta_1 + \epsilon
\ \land \ \cdots \ \land \ \delta_k \lt \tilde\lambda_k-\tilde\lambda_{k-1} \lt \delta_k + \epsilon)\text{,}\label{sec_stable-2-8}
\end{equation}
where \(e^{i \lambda_0},e^{i\lambda_1},\ldots,e^{i\lambda_k}\) are consecutive eigenvalues of a random matrix in~\(U(N)\) and \(\tilde\lambda = N\lambda/2\pi\) are the normalized eigenangles. In other words, the value of the PDF of the normalized \(k\)-neighbor gaps. By inspection of Haar measure in terms of eigenvalues, the likelihood of a configuration is nonzero if \(N\) is sufficiently large, and as \(N\to\infty\) that likelihood converges to a nonzero limit. Furthermore, that likelihood can be written in terms of the \(m\)-correlation functions, with \(m \le k\), and each \(m\)-correlation converges to its limit at the rate of \(1/N^2\).  Therefore the likelihood also converges at that rate.%
\par
Thus, if a property of a unitary polynomial depends only on the local configuration of zeros, then the limiting behavior of that property can be reliably studied by performing numerical experiments on random matrices from \(U(N)\) for moderately sized~\(N\). Properties which rely on large numbers of zeros, such as moments or value distribution, can show significant departures from the limiting behavior.%
%
%
\typeout{************************************************}
\typeout{Subsection 1.4 Curvature effects}
\typeout{************************************************}
\subsection{Curvature effects}\label{sec_curve}
As \(N\to\infty\) a given configuration sits on a segment of the circle of length \(O(1/N)\).  That segment is approaching a straight line, and the discrepancy is~\(O(1/N^2)\).  So that effect is comparable to the previously discussed effects.%
%
%
\typeout{************************************************}
\typeout{Section 2 The source of the second bump}
\typeout{************************************************}
\section{The source of the second bump}\label{thesource}
In this section we justify the following:%
\begin{claim}\label{sourceofhump2}%
The zeros of \(p'\) in the second bump, which do not come from two closely spaced zeros, primarily arise from three consecutive zeros where both zero gaps are larger than average.  The associated zero of the derivative is near the radial line from the center of the circle to the middle zero.%
\par
Specifically, if \(e^{i \theta_1}, e^{i \theta_2}, e^{i \theta_3}\) are consecutive zeros of~\(p\), with \(\theta_2 - \theta_1, \theta_3 - \theta_2 > 2\pi/N\), then \(p'\) tends to have a zero within \(2.5/N\) of~\((1-3/N) z_2\).%
\end{claim}
\hyperref[specialcircles]{Figure~{\xreffont\ref{specialcircles}}} illustrates \hyperref[sourceofhump2]{Claim~{\xreffont\ref{sourceofhump2}}} using an example with \(N=30\).%
\par
That claim reinforces the idea that the distribution of zeros of \(p'\) is difficult to analyze.  Given three zeros of \(p\), we have located one zero of \(p'\). That leaves at least one, and typically two, zeros unaccounted for. (But when two consecutive gaps are larger than average, it is more likely that the gaps on either side are smaller than average, so that says something about the location of the ``missing'' zeros.  \hyperref[specialcircles]{Figure~{\xreffont\ref{specialcircles}}} contains some examples.)%
\begin{figureptx}{Figure}{The zeros of a degree \(N=30\) unitary polynomial with circles of radius \(2.5/N\) containing some of the zeros of the derivative.}{specialcircles}{}%
\begin{image}{0.175}{0.65}{0.175}{}%
\includegraphics[width=\linewidth]{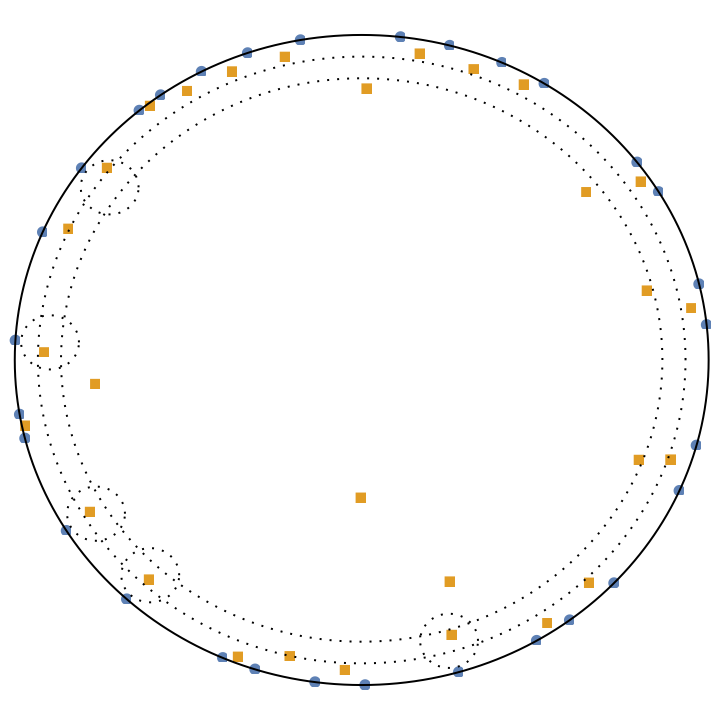}
\end{image}%
\tcblower
\end{figureptx}%
We reinforce \hyperref[sourceofhump2]{Claim~{\xreffont\ref{sourceofhump2}}} with the following:%
\begin{example}[Focusing on the special zeros.]\label{exp_special20}%
Start with an empty list \mono{specialzeros20}. Generate approximately \(100,000\) characteristic polynomials \(p_u\) for Haar random \(u\in U(20)\).  In all cases where three consecutive zeros \(z_1, z_2, z_3\) have both gaps larger than average, select all zeros of the derivative in the target disc of radius \(2.5/N\) centered at \((1-3/N) z_2\). (That is, select the \(z'\) in the small circles, as in \hyperref[specialcircles]{Figure~{\xreffont\ref{specialcircles}}}.) Append those zeros to \mono{specialzeros20}.%
\end{example}
In \hyperref[exp_special20]{Experiment~{\xreffont\ref{exp_special20}}}, the condition on consecutive zeros occurred around \(332,000\) times, which is consistent with the expected frequency of consecutive gaps which are larger than average. In those cases, in the target disc, \(99.43\)\% of the time there was exactly one zero. Of the remaining \(0.57\)\%, \(92.7\)\% of the time there were two zeros, and in the remaining \(0.04\)\%, no zeros.%
\par
\hyperref[pointsinspecialcircles]{Figure~{\xreffont\ref{pointsinspecialcircles}}} shows all the zeros from \hyperref[exp_special20]{Experiment~{\xreffont\ref{exp_special20}}} in the target discs, as individual points and as a density plot. That figure gives a clear picture that the configuration of consecutive larger-than-average gaps contributes significantly to the second maximum.%
\begin{figureptx}{Figure}{The points from \mono{specialzeros20} in \hyperref[exp_special20]{Experiment~{\xreffont\ref{exp_special20}}}, with the target disc rotated so that the corresponding zero \(z_2\) lies at the black dot. Point cloud on the left and histogram on the right.}{pointsinspecialcircles}{}%
\begin{sidebyside}{2}{0}{0}{0.06}%
\begin{sbspanel}{0.47}%
\includegraphics[width=\linewidth]{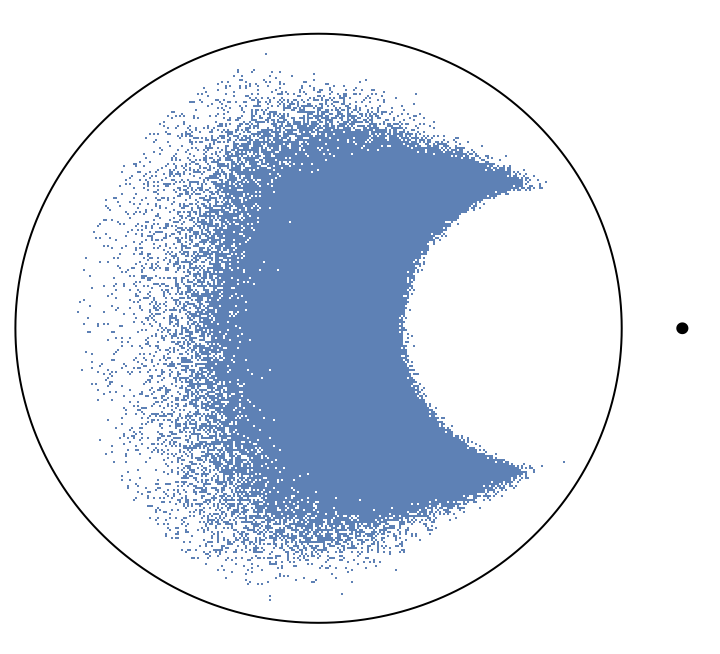}
\end{sbspanel}%
\begin{sbspanel}{0.47}%
\includegraphics[width=\linewidth]{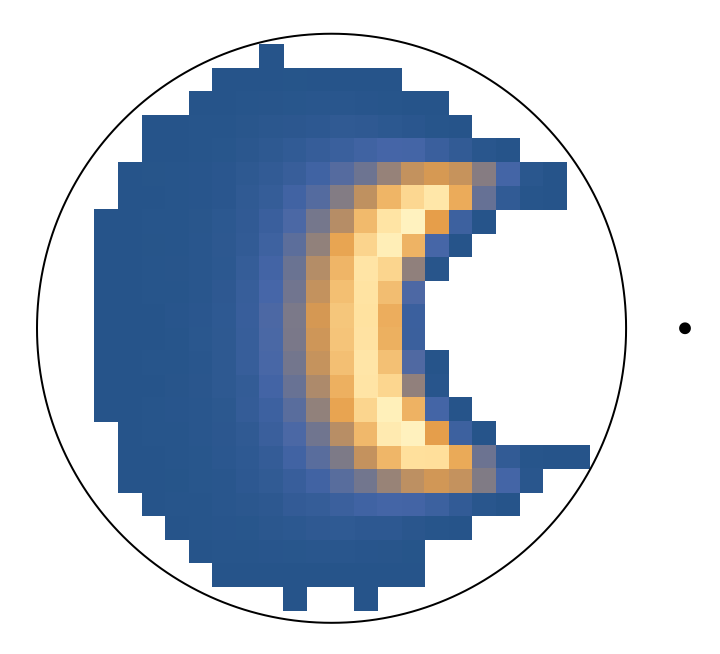}
\end{sbspanel}%
\end{sidebyside}%
\tcblower
\end{figureptx}%
\hyperref[specialhistogram]{Figure~{\xreffont\ref{specialhistogram}}} shows the histogram, normalized to have area~\(1\), of \(20(1-\abs{z'})\) for the points \(z'\) in \hyperref[pointsinspecialcircles]{Figure~{\xreffont\ref{pointsinspecialcircles}}}. Considering the arrangement of the points in \hyperref[pointsinspecialcircles]{Figure~{\xreffont\ref{pointsinspecialcircles}}}, one should not expect that density function to have a simple form. Note that the heavier tail to the right in the histogram corresponds to the left portion of the cloud of dots.%
\begin{figureptx}{Figure}{Histogram of \(20(1-\abs{z'})\) for the \(z'\) in \hyperref[pointsinspecialcircles]{Figure~{\xreffont\ref{pointsinspecialcircles}}}}{specialhistogram}{}%
\begin{image}{0.25}{0.5}{0.25}{}%
\includegraphics[width=\linewidth]{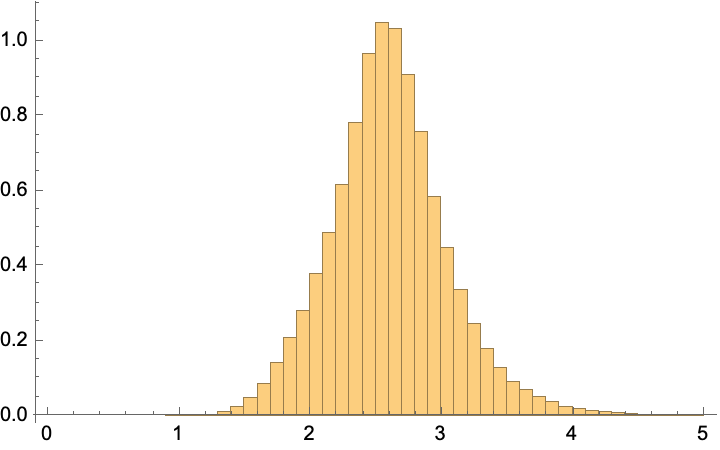}
\end{image}%
\tcblower
\end{figureptx}%
%
%
\typeout{************************************************}
\typeout{Section 3 A toy model}
\typeout{************************************************}
\section{A toy model}\label{toymodel}
We have argued that the second bump arises when consecutive zero gaps are larger than average, with the details about the remaining zeros playing a lesser role.  We will explore this by analyzing a simple example.%
%
%
\typeout{************************************************}
\typeout{Subsection 3.1 Excised roots of unity}
\typeout{************************************************}
\subsection{Excised roots of unity}\label{toymodel-3}
Consider \(x^n-1\) where two zeros have been deleted, creating consecutive zero gaps of twice the average.  That is, let%
\begin{equation}
f_n(x) = \frac{x^n-1}{(1-e^{2\pi i/n}x)(1-e^{-2\pi i/n}x)}
= \frac{x^n-1}{1-2 \cos(2 \pi/n)x + x^2}\text{.}\label{toymodel-3-2-2}
\end{equation}
\hyperref[model1fig]{Figure~{\xreffont\ref{model1fig}}} illustrates the cases \(n=20\) and \(n=40\).%
\begin{figureptx}{Figure}{The zeros, and the zeros of the derivative, of \(f_{20}\) on the left and \(f_{40}\) on the right.}{model1fig}{}%
\begin{sidebyside}{2}{0}{0}{0.1}%
\begin{sbspanel}{0.45}%
\includegraphics[width=\linewidth]{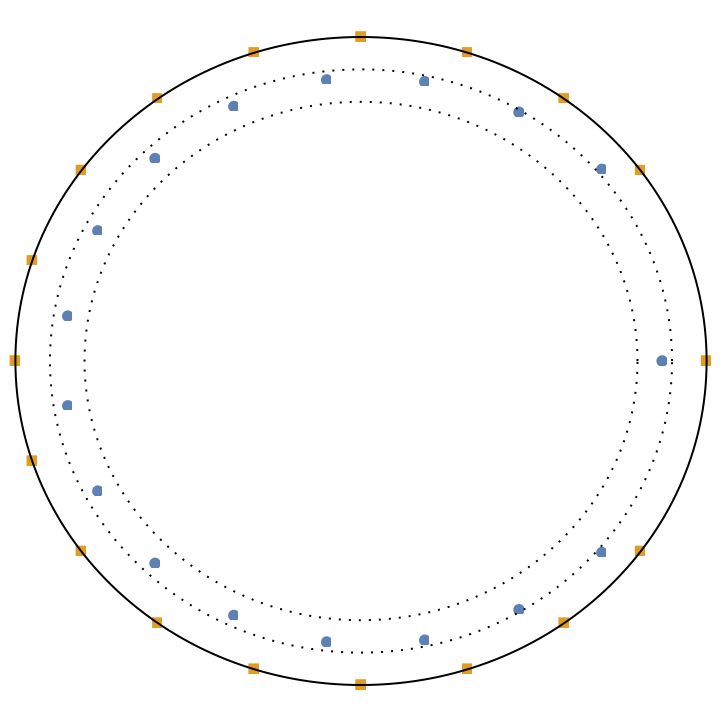}
\end{sbspanel}%
\begin{sbspanel}{0.45}%
\includegraphics[width=\linewidth]{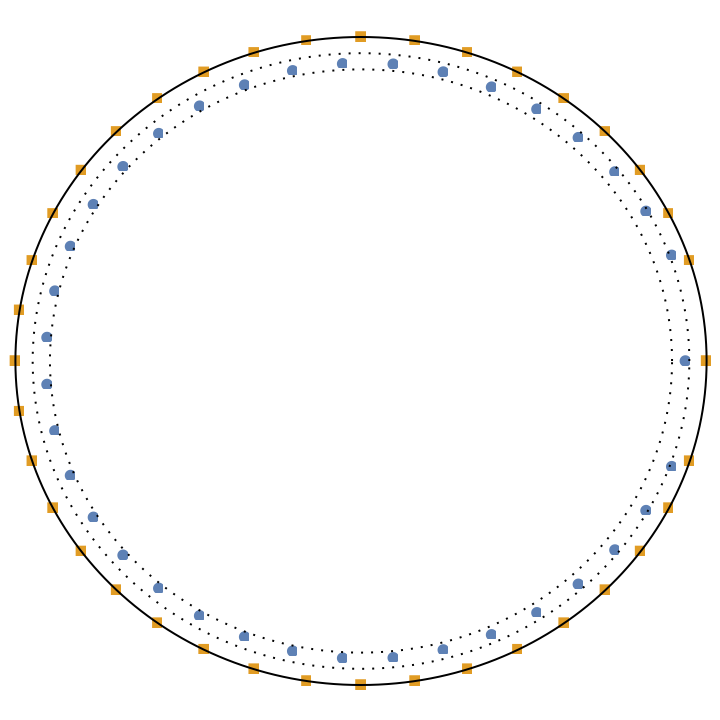}
\end{sbspanel}%
\end{sidebyside}%
\tcblower
\end{figureptx}%
\begin{proposition}\label{toymodel-3-4}%
The derivative of \(f_n\) has a zero within \(O(1/n^2)\) of \(1 - b_0/n\), where \(b_0 = 2.3565\ldots\) is the unique real root of \(4 \pi^2 +2 b + b^2 - 2 b e^b\).%
\end{proposition}
\begin{proof}\label{toymodel-3-4-2}
By the quotient rule, the zeros of \(f'_n\) are the same as the zeros of%
\begin{equation}
F_n(x) = -2 \cos(2 \pi/n)  + 2 x + n x^{-1 + n} + 2 \cos(2 \pi/n)(1-n)  x^n - (2-n) x^{1 + n}\label{toymodel-3-4-2-1-2}
\end{equation}
Since \(F_n(0) \approx -2\) and \(F_n(1) \approx 4\pi^2/n\), we see that \(F_n\) has a real zero in the interval \((0, 1)\), and that zero should be close to \(1\).%
\par
Setting \(x = 1-b/n\) and expanding to order \(1/n^2\) we find%
\begin{equation}
F_n(1 - b/n) =
e^{-b}\, \frac{2 b + b^2 - 2 b e^b + 4 \pi^2}n + 
e^{-b} \, \frac{-b^4 - 8 \pi^2 - 4 b^2 \pi^2 + 8 e^b \pi^2}{2 n^2}
+O\left(\frac{1}{n^3}\right) \text{.}\label{toymodel-3-4-2-2-3}
\end{equation}
Thus, \(F_n(1 - b_0/n) = O(1/n^2)\). The argument is completed by Rouché's theorem.%
\end{proof}
Suppose the polynomial \(f_n\) was modified so that the roots at \(e^{\pm 2\pi i \cdot 2/n}\) instead were at \(e^{2\pi i \cdot c_\pm /n}\), where \(1\lt c_+ \lt 3\) and \(-3\lt c_- \lt -1\). Using the above method, one can check throughout those ranges the derivative of the modified polynomial  has a root within \(0.8/n\) of \(1-2.6/n\). This is a small subset of the region asserted in \hyperref[sourceofhump2]{Claim~{\xreffont\ref{sourceofhump2}}}.%
%
%
\typeout{************************************************}
\typeout{Subsection 3.2 A final remark}
\typeout{************************************************}
\subsection{A final remark}\label{toymodel-4}
We have seen that, from the perspective of the geometry of zeros of polynomials, the behavior of the zeros of the derivative of a unitary polynomial \(p\) is quite complicated.  There are at least three (partially overlapping) regimes for the zeros of~\(p'\), and in each regime the nature of the dominant influence is different.%
\par
This suggests that making progress on understanding the distribution of zeros of~\(p'\), such as proving that the distribution is bimodal, may need to come from another approach.  The expected value of moments of \(\abs{p_u'(1-r/N)}\) for \(u\in U(N)\) seems a more likely avenue for progress. However, it could still be interesting to quantify the way in which appropriately spaced triples of zeros contribute to the 2nd mode.%
\subparagraph*{Acknowledgments.}\label{toymodel-5-1}%
This work was supported by the National Science Foundation. This paper was written in PreTeXt~\hyperlink{PTX}{[{\xreffont 5}]}.%
\subparagraph*{Disclosure.}\label{toymodel-5-2}%
The author reports there are no competing interests to declare.%
\addtocontents{toc}{\vspace{\normalbaselineskip}}
\typeout{************************************************}
\typeout{References  Bibliography}
\typeout{************************************************}
\bibliographystyle{amsplain}

\end{document}